\documentclass[12pt]{article}
\usepackage{amssymb}
\usepackage{amsfonts}
\usepackage{amsmath}
\usepackage{amsbsy}
\usepackage[small,nohug,heads=littlevee]{diagrams}
\diagramstyle[labelstyle=\scriptstyle]
\newcommand{\be}{\begin{equation}}
\newcommand{\ee}{\end{equation}}
\newcommand{\bea}{\begin{eqnarray}}
\newcommand{\eea}{\end{eqnarray}}
\newcommand{\ba}{\begin{array}}
\newcommand{\ea}{\end{array}}

\newcommand{\bc}{\begin{center}}
\newcommand{\ec}{\end{center}}
\newcommand{\ben}{\begin{enumerate}}
\newcommand{\een}{\end{enumerate}}
\newcommand{\bfi}{\begin{figure}}
\newcommand{\efi}{\end{figure}}

\newcommand{\bq}{\begin{quote}}
\newcommand{\eq}{\end{quote}}
\newcommand{\bqu}{\begin{quotation}}
\newcommand{\equ}{\end{quotation}}
\newenvironment{emphit}{\begin{itemize}}{\end{itemize}}
\newcommand{\bemp}{\begin{emphit}}
\newcommand{\eemp}{\end{emphit}}

\newcommand{\bt}{\begin{tabular}}
\newcommand{\et}{\end{tabular}}

\newtheorem{myth}{Theorem}[section]
\newtheorem{mylem}{Lemma}[section]
\newtheorem{mycor}{Corollary}[section]

\newtheorem{mydef}{Definition}[section]
\newtheorem{myexam}{Example}[section]
\newtheorem{myquest}{Question}[section]
\newtheorem{myrem}{Remark}[section]

\begin{document}
\date{}
\title{Smarandache Isotopy Theory Of Smarandache: Quasigroups And Loops
\footnote{2000 Mathematics Subject Classification. Primary 20NO5 ;
Secondary 08A05.}
\thanks{{\bf Keywords and Phrases : } Smarandache: groupoids; quasigroups; loops; $f$,$g$ principal isotopes}}
\author{T\`em\'it\d{\'o}p\d{\'e} Gb\d{\'o}l\'ah\`an Ja\'iy\'e\d ol\'a\\Department of
Mathematics,\\
Obafemi Awolowo University, Ile Ife,
Nigeria.\\jaiyeolatemitope@yahoo.com, tjayeola@oauife.edu.ng}
\maketitle

\begin{abstract}
The concept of Smarandache isotopy is introduced and its study is
explored for Smarandache: groupoids, quasigroups and loops just like
the study of isotopy theory was carried out for groupoids,
quasigroups and loops. The exploration includes: Smarandache;
isotopy and isomorphy classes, Smarandache $f,g$ principal isotopes
and G-Smarandache loops.
\end{abstract}

\section{Introduction}
\paragraph{}In 2002, W. B. Vasantha Kandasamy initiated the study of
Smarandache loops in her book \cite{phd75} where she introduced over
75 Smarandache concepts in loops. In her paper \cite{phd83}, she
defined a Smarandache loop (S-loop) as a loop with at least a
subloop which forms a subgroup under the binary operation of the
loop. For more on loops and their properties, readers should check
\cite{phd3}, \cite{phd41}, \cite{phd39}, \cite{phd49}, \cite{phd42}
and \cite{phd75}. In [\cite{phd75}, Page~102], the author introduced
Smarandache isotopes of loops particularly Smarandache principal
isotopes. She has also introduced the Smarandache concept in some other algebraic structures
as \cite{van1,van2,van3,van4,van5,van6} account. The present author has contributed to the study of
S-quasigroups and S-loops in \cite{sma1}, \cite{sma2} and
\cite{sma3} while Muktibodh \cite{muk} did a study on the first.

In this study, the concept of Smarandache isotopy will be introduced
and its study will be explored in Smarandache: groupoids,
quasigroups and loops just like the study of isotopy theory was
carried out for groupoids, quasigroups and loops as summarized in
Bruck \cite{phd41}, Dene and Keedwell \cite{phd49}, Pflugfelder
\cite{phd3}.

\section{Definitions and Notations}
\begin{mydef}\label{1:0}
Let $L$ be a non-empty set. Define a binary operation ($\cdot $) on
$L$ : If $x\cdot y\in L~\forall ~x, y\in L$, $(L, \cdot )$ is called
a groupoid. If the system of equations ; $a\cdot x=b$ and $y\cdot
a=b$ have unique solutions for $x$ and $y$ respectively, then $(L,
\cdot )$ is called a quasigroup. Furthermore, if there exists a
unique element $e\in L$ called the identity element such that
$\forall
~x\in L$, $x\cdot e=e\cdot x=x$, $(L, \cdot )$ is called a loop.\\

If there exists at least a non-empty and non-trivial subset $M$ of a
groupoid(quasigroup or semigroup or loop) $L$ such that $(M,\cdot )$
is a non-trivial subsemigroup(subgroup or subgroup or subgroup) of
$(L, \cdot )$, then $L$ is called a Smarandache:
groupoid(S-groupoid)$\Big($quasigroup(S-quasigroup) or
semigroup(S-semigroup) or loop(S-loop)$\Big)$ with Smarandache:
subsemigroup(S-subsemigroup)$\Big($subgroup(S-subgroup) or
subgroup(S-subgroup) or subgroup(S-subgroup)$\Big)$ $M$.

Let $(G,\cdot )$ be a quasigroup(loop). The bijection $L_x :
G\rightarrow G$ defined as $yL_x=x\cdot y~\forall ~x, y\in G$ is
called a left translation(multiplication) of $G$ while the bijection
$R_x : G\rightarrow G$ defined as $yR_x=y\cdot x~\forall ~x, y\in G$
is called a right translation(multiplication) of $G$.

The set $SYM(L, \cdot )=SYM(L)$ of all bijections in a groupoid
$(L,\cdot )$ forms a group called the permutation(symmetric) group
of the groupoid $(L,\cdot )$.
\end{mydef}

\begin{mydef}\label{1:1}
If $(L, \cdot )$ and $(G, \circ )$ are two distinct groupoids, then
the triple $(U, V, W) : (L, \cdot )\rightarrow (G, \circ )$ such
that $U, V, W : L\rightarrow G$ are bijections is called an
isotopism if and only if
\begin{displaymath}
xU\circ yV=(x\cdot y)W~\forall ~x, y\in L.
\end{displaymath}
So we call $L$ and $G$ groupoid isotopes. If $L=G$ and
$W=I$(identity mapping) then $(U, V, I)$ is called a principal
isotopism, so we call $G$ a principal isotope of $L$. But if in
addition $G$ is a quasigroup such that for some $f,g\in G$, $U=R_g$
and $V=L_f$, then $(R_g,L_f,I)~:~(G, \cdot )\rightarrow (G, \circ )$
is called an $
f,g$-principal isotopism while $(G, \cdot )$ and $(G,
\circ )$ are called quasigroup isotopes.

If $U=V=W$, then $U$ is called an isomorphism, hence we write $(L,
\cdot )\cong (G, \circ )$. A loop $(L, \cdot )$ is called a G-loop
if and only if $(L, \cdot )\cong (G, \circ )$ for all loop isotopes
$(G, \circ )$ of $(L, \cdot )$.

Now, if $(L, \cdot )$ and $(G, \circ )$ are S-groupoids with
S-subsemigroups $L'$ and $G'$ respectively such that $(G')A=L'$, where $A\in\{U,V,W\}$, then the isotopism $(U,
V, W) : (L, \cdot )\rightarrow (G, \circ )$ is called a Smarandache
isotopism(S-isotopism). Consequently, if $W=I$ the triple $(U, V,
I)$ is called a Smarandache principal isotopism. But if in addition
$G$ is a S-quasigroup with S-subgroup $H'$ such that for some
$f,g\in H$, $U=R_g$ and $V=L_f$, and $(R_g,L_f,I)~:~(G, \cdot
)\rightarrow (G, \circ )$ is an isotopism, then the triple is called
a Smarandache $f,g$-principal isotopism while $f$ and $g$ are called
Smarandache elements(S-elements).

Thus, if $U=V=W$, then $U$ is called a Smarandache isomorphism,
hence we write $(L, \cdot )\succsim (G, \circ )$. An S-loop $(L,
\cdot )$ is called a G-Smarandache loop(GS-loop) if and only if $(L,
\cdot )\succsim (G, \circ )$ for all loop isotopes(or particularly
all S-loop isotopes) $(G, \circ )$ of $(L, \cdot )$.
\end{mydef}

\begin{myexam}\label{exam:1}
The systems $(L,\cdot )$ and $(L,\ast)$, $L=\{0,1,2,3,4\}$ with the multiplication tables below are S-quasigroups with
S-subgroups $(L',\cdot )$ and $(L'',\ast )$ respectively, $L'=\{0,1\}$ and $L''=\{1,2\}$. $(L,\cdot )$ is taken from Example~2.2 of \cite{muk}.
The triple $(U,V,W)$ such that
\begin{displaymath}
U=\left(\begin{array}{ccccc} 0 & 1 & 2 & 3 & 4 \\
1 & 2 & 3 & 4 & 0
\end{array}\right),~V=
\left(\begin{array}{ccccc}
0 & 1 & 2 & 3 & 4 \\
1 & 2 & 4 & 0 & 3
\end{array}\right)~\textrm{and}~W=
\left(\begin{array}{ccccc}
0 & 1 & 2 & 3 & 4 \\
1 & 2 & 0 & 4 & 3
\end{array}\right)
\end{displaymath}
are permutations on $L$, is an S-isotopism of $(L,\cdot )$ onto $(L,\ast)$. Notice that $A(L')=L''$ for all $A\in\{U,V,W\}$
and $U, V, W : L'\rightarrow L''$ are all bijcetions.\\
\begin{center}
\begin{tabular}{|c||c|c|c|c|c|}
\hline
$\cdot $ & 0 & 1 & 2 & 3 & 4 \\
\hline \hline
0 & 0 & 1 & 3 & 4 & 2 \\
\hline
1 & 1 & 0 & 2 & 3 & 4 \\
\hline
2 & 3 & 4 & 1 & 2 & 0 \\
\hline
3 & 4 & 2 & 0 & 1 & 3 \\
\hline
4 & 2 & 3 & 4 & 0 & 1 \\
\hline
\end{tabular}
\qquad\qquad\qquad\qquad\qquad\qquad\qquad\qquad\begin{tabular}{|c||c|c|c|c|c|}
\hline
$\ast$ & 0 & 1 & 2 & 3 & 4 \\
\hline \hline
0 & 1 & 0 & 4 & 2 & 3 \\
\hline
1 & 3 & 1 & 2 & 0 & 4 \\
\hline
2 & 4 & 2 & 1 & 3 & 0 \\
\hline
3 & 0 & 4 & 3 & 1 & 2 \\
\hline
4 & 2 & 3 & 0 & 4 & 1 \\
\hline
\end{tabular}
\end{center}
\end{myexam}

\begin{myexam}\label{exam:2}
According to Example~4.2.2 of \cite{van2}, the system $(\mathbb{Z}_6,\times_6)$ i.e the set $L=\mathbb{Z}_6$ under multiplication modulo $6$ is an S-semigroup
with S-subgroups $(L',\times_6 )$ and $(L'',\times_6 )$, $L'=\{2,4\}$ and $L''=\{1,5\}$. This can be deduced from its multiplication table, below.
The triple $(U,V,W)$ such that
\begin{displaymath}
U=\left(\begin{array}{cccccc} 0 & 1 & 2 & 3 & 4 & 5\\
4 & 3 & 5 & 1 & 2 & 0
\end{array}\right),~V=
\left(\begin{array}{cccccc}
0 & 1 & 2 & 3 & 4 & 5 \\
1 & 3 & 2 & 4 & 5 & 0
\end{array}\right)~\textrm{and}~W=
\left(\begin{array}{cccccc}
0 & 1 & 2 & 3 & 4 & 5\\
1 & 0 & 5 & 4 & 2 & 3
\end{array}\right)
\end{displaymath}
are permutations on $L$, is an S-isotopism of $(\mathbb{Z}_6,\times_6)$ unto an S-semigroup $(\mathbb{Z}_6,\ast)$ with
S-subgroups $(L''',\ast )$ and $(L'''',\ast )$, $L'''=\{2,5\}$ and $L''''=\{0,3\}$ as shown in the second table below.
Notice that $A(L')=L'''$ and $A(L'')=L''''$ for all $A\in\{U,V,W\}$
and $U, V, W : L'\rightarrow L'''$ and $U, V, W : L''\rightarrow L''''$ are all bijcetions.\\
\begin{center}
\begin{tabular}{|c||c|c|c|c|c|c|}
\hline
$\times_6$ & 0 & 1 & 2 & 3 & 4 & 5\\
\hline \hline
0 & 0 & 0 & 0 & 0 & 0 & 0 \\
\hline
1 & 0 & 1 & 2 & 3 & 4 & 5 \\
\hline
2 & 0 & 2 & 4 & 0 & 2 & 4\\
\hline
3 & 0 & 3 & 0 & 3 & 0 & 3 \\
\hline
4 & 0 & 4 & 2 & 0 & 4 & 2 \\
\hline
5 & 0 & 5 & 4 & 3 & 2 & 1 \\
\hline
\end{tabular}
\qquad\qquad\qquad\qquad\qquad\qquad\qquad\qquad\begin{tabular}{|c||c|c|c|c|c|c|}
\hline
$\ast$ & 0 & 1 & 2 & 3 & 4 & 5\\
\hline \hline
0 & 0 & 1 & 2 & 3 & 4 & 5 \\
\hline
1 & 4 & 1 & 1 & 4 & 4 & 1 \\
\hline
2 & 5 & 1 & 5 & 2 & 1 & 2\\
\hline
3 & 3 & 1 & 5 & 0 & 4 & 2 \\
\hline
4 & 1 & 1 & 1 & 1 & 1 & 1 \\
\hline
5 & 2 & 1 & 2 & 5 & 1 & 5 \\
\hline
\end{tabular}
\end{center}
\end{myexam}

\begin{myrem}
Taking careful look at Definition~\ref{1:1} and comparing it with
[Definition~4.4.1,\cite{phd75}], it will be observed that the author
did not allow the component bijections $U$,$V$ and $W$ in $(U, V,
W)$ to act on the whole S-loop $L$ but only on the
S-subloop(S-subgroup) $L'$. We feel this is necessary to adjust here
so that the set $L-L'$ is not out of the study. Apart from this, our
adjustment here will allow the study of Smarandache isotopy to be
explorable. Therefore, the S-isotopism and S-isomorphism here are
clearly special types of relations(isotopism and isomorphism) on the
whole domain into the whole co-domain but those of Vasantha
Kandasamy \cite{phd75} only take care of the structure of the
elements in the S-subloop and not the S-loop. Nevertheless, we do
not fault her study for we think she defined them to apply them to
some life problems as an applied algebraist.
\end{myrem}

\section{Smarandache Isotopy and Isomorphy Classes}
\begin{myth}\label{1:3}
Let $\mathfrak{G}=\Big\{\big(G_\omega,\circ_\omega
\big)\Big\}_{\omega\in \Omega}$ be a set of distinct S-groupoids
with a corresponding set of S-subsemigroups
$\mathfrak{H}=\Big\{\big(H_\omega,\circ_\omega
\big)\Big\}_{\omega\in \Omega}$. Define a relation $\thicksim$ on
$\mathfrak{G}$ such that for all $\big(G_{\omega_i},\circ_{\omega_i}
\big)~,\big(G_{\omega_j},\circ_{\omega_j} \big)\in\mathfrak{G}$,
where $\omega_i,\omega_j\in \Omega$,
\begin{displaymath}
\big(G_{\omega_i},\circ_{\omega_i} \big)\thicksim
\big(G_{\omega_j},\circ_{\omega_j} \big)\Longleftrightarrow
\big(G_{\omega_i},\circ_{\omega_i}
\big)~\textrm{and}~\big(G_{\omega_j},\circ_{\omega_j}
\big)~\textrm{are S-isotopic}.
\end{displaymath}
Then $\thicksim$ is an equivalence relation on $\mathfrak{G}$.
\end{myth}
{\bf Proof}\\
Let $\big(G_{\omega_i},\circ_{\omega_i}
\big)~,\big(G_{\omega_j},\circ_{\omega_j}
\big),~\big(G_{\omega_k},\circ_{\omega_k} \big),\in\mathfrak{G}$,
where $\omega_i,\omega_j,\omega_k\in \Omega$.
\begin{description}
\item[Reflexivity] If $I~:~G_{\omega_i}\rightarrow G_{\omega_i}$ is the
identity mapping, then
\begin{displaymath}
xI\circ_{\omega_i}yI=(x\circ_{\omega_i} y)I~\forall~x,y\in
G_{\omega_i}\Longrightarrow~\textrm{the triple}~
(I,I,I)~:~\big(G_{\omega_i},\circ_{\omega_i} \big)\rightarrow
\big(G_{\omega_i},\circ_{\omega_i} \big)
\end{displaymath}
is an S-isotopism since
$\big(H_{\omega_i}\big)I=H_{\omega_i}~\forall~\omega_i\in \Omega$.
In fact, it can be simply deduced that every S-groupoid is
S-isomorphic to itself.
\item[Symmetry] Let $\big(G_{\omega_i},\circ_{\omega_i} \big)\thicksim
\big(G_{\omega_j},\circ_{\omega_j} \big)$. Then there exist
bijections
\begin{displaymath}
U,V,W~:~\big(G_{\omega_i},\circ_{\omega_i} \big)\longrightarrow
\big(G_{\omega_j},\circ_{\omega_j} \big)~\textrm{such that}~
\big(H_{\omega_i}\big)A=H_{\omega_j}~\forall~A\in\{U,V,W\}
\end{displaymath} so that
the triple
\begin{displaymath}
\alpha =(U, V, W)~:~\big(G_{\omega_i},\circ_{\omega_i}
\big)\longrightarrow \big(G_{\omega_j},\circ_{\omega_j} \big)
\end{displaymath}
is an isotopism. Since each of $U,V,W$ is bijective, then their
inverses
\begin{displaymath}
U^{-1},V^{-1},W^{-1}~:~\big(G_{\omega_j},\circ_{\omega_j}
\big)\longrightarrow \big(G_{\omega_i},\circ_{\omega_i} \big)
\end{displaymath} are
bijective. In fact,
$\big(H_{\omega_j}\big)A^{-1}=H_{\omega_i}~\forall~A\in\{U,V,W\}$
since $A$ is bijective so that the triple
\begin{displaymath}
\alpha^{-1} =(U^{-1}, V^{-1},
W^{-1})~:~\big(G_{\omega_j},\circ_{\omega_j} \big)\longrightarrow
\big(G_{\omega_i},\circ_{\omega_i} \big)
\end{displaymath}
is an isotopism. Thus, $\big(G_{\omega_j},\circ_{\omega_j}
\big)\thicksim \big(G_{\omega_i},\circ_{\omega_i} \big)$.
\item[Transitivity] Let $\big(G_{\omega_i},\circ_{\omega_i} \big)\thicksim
\big(G_{\omega_j},\circ_{\omega_j} \big)$ and
$\big(G_{\omega_j},\circ_{\omega_j} \big)\thicksim
\big(G_{\omega_k},\circ_{\omega_k} \big)$. Then there exist
bijections
\begin{displaymath}
U_1,V_1,W_1~:~\big(G_{\omega_i},\circ_{\omega_i}
\big)\longrightarrow \big(G_{\omega_j},\circ_{\omega_j}
\big)~\textrm{and}~U_2,V_2,W_2~:~\big(G_{\omega_j},\circ_{\omega_j}
\big)\longrightarrow \big(G_{\omega_k},\circ_{\omega_k} \big)
\end{displaymath}
\begin{displaymath}
\textrm{such that}~
\big(H_{\omega_i}\big)A=H_{\omega_j}~\forall~A\in\{U_1,V_1,W_1\}
\end{displaymath}
\begin{displaymath}
\textrm{and}~\big(H_{\omega_j}\big)B=H_{\omega_k}~\forall~B\in\{U_2,V_2,W_2\}~\textrm{so
that the triples}
\end{displaymath}
\begin{displaymath}
\alpha_1 =(U_1,V_1,W_1)~:~\big(G_{\omega_i},\circ_{\omega_i}
\big)\longrightarrow \big(G_{\omega_j},\circ_{\omega_j}
\big)~\textrm{and}
\end{displaymath}
\begin{displaymath}
\alpha_2 =(U_2,V_2,W_2)~:~\big(G_{\omega_j},\circ_{\omega_j}
\big)\longrightarrow \big(G_{\omega_k},\circ_{\omega_k} \big)
\end{displaymath} are
isotopisms. Since each of $U_i,V_i,W_i$, $i=1,2$, is bijective, then
\begin{displaymath}
U_3=U_1U_2,V_3=V_1V_2,W_3=W_1W_2~:~\big(G_{\omega_i},\circ_{\omega_i}
\big)\longrightarrow \big(G_{\omega_k},\circ_{\omega_k} \big)
\end{displaymath}
are bijections such that
$\big(H_{\omega_i}\big)A_3=\big(H_{\omega_i}\big)A_1A_2=\big(H_{\omega_j}\big)A_2=H_{\omega_k}$
so that the triple
\begin{displaymath}
\alpha_3=\alpha_1\alpha_2=(U_3,V_3,W_3)~:~\big(G_{\omega_i},\circ_{\omega_i}
\big)\longrightarrow \big(G_{\omega_k},\circ_{\omega_k} \big)
\end{displaymath}
is an isotopism. Thus, $\big(G_{\omega_i},\circ_{\omega_i}
\big)\thicksim \big(G_{\omega_k},\circ_{\omega_k}\big)$.
\end{description}

\begin{myrem}
As a follow up to Theorem~\ref{1:3}, the elements of the set
$\mathfrak{G}/\thicksim$ will be referred to as Smarandache isotopy
classes(S-isotopy classes). Similarly, if $\thicksim$ meant
"S-isomorphism" in Theorem~\ref{1:3}, then the elements of
$\mathfrak{G}/\thicksim$ will be referred to as Smarandache
isomorphy classes(S-isomorphy classes). Just like isotopy has an
advantage over isomorphy in the classification of loops, so also
S-isotopy will have advantage over S-isomorphy in the classification
of S-loops.
\end{myrem}

\begin{mycor}\label{1:4}
Let $\mathcal{L}_n$, $\mathcal{SL}_n$ and $\mathcal{NSL}_n$ be the
sets of; all finite loops of order $n$; all finite S-loops of order
$n$ and all finite non S-loops of order $n$ respectively.
\begin{enumerate}
\item If
$\mathcal{A}_i^n$ and $\mathcal{B}_i^n$ represent the isomorphy
class of $\mathcal{L}_n$ and the S-isomorphy class of
$\mathcal{SL}_n$ respectively, then
\begin{description}
\item[(a)] $|\mathcal{SL}_n|+|\mathcal{NSL}_n|=|\mathcal{L}_n|$;
\begin{description}
\item[(i)] $|\mathcal{SL}_5|+|\mathcal{NSL}_5|=56$,
\item[(ii)] $|\mathcal{SL}_6|+|\mathcal{NSL}_6|=9, 408$ and
\item[(iii)] $|\mathcal{SL}_7|+|\mathcal{NSL}_7|=16, 942,080$.
\end{description}
\item[(b)] $|\mathcal{NSL}_n|=\sum_{i=1}|\mathcal{A}_i^n|-\sum_{i=1}|\mathcal{B}_i^n|$;
\begin{description}
\item[(i)]
$|\mathcal{NSL}_5|=\sum_{i=1}^6|\mathcal{A}_i^5|-\sum_{i=1}|\mathcal{B}_i^5|$,
\item[(ii)]
$|\mathcal{NSL}_6|=\sum_{i=1}^{109}|\mathcal{A}_i^6|-\sum_{i=1}|\mathcal{B}_i^6|$
and
\item[(iii)] $|\mathcal{NSL}_7|=\sum_{i=1}^{23,
746}|\mathcal{A}_i^7|-\sum_{i=1}|\mathcal{B}_i^7|$.
\end{description}
\end{description}
\item If
$\mathfrak{A}_i^n$ and $\mathfrak{B}_i^n$ represent the isotopy
class of $\mathcal{L}_n$ and the S-isotopy class of $\mathcal{SL}_n$
respectively, then
\begin{displaymath}
|\mathcal{NSL}_n|=\sum_{i=1}|\mathfrak{A}_i^n|-\sum_{i=1}|\mathfrak{B}_i^n|;
\end{displaymath}
\begin{description}
\item[(i)]
$|\mathcal{NSL}_5|=\sum_{i=1}^2|\mathfrak{A}_i^5|-\sum_{i=1}|\mathfrak{B}_i^5|$,
\item[(ii)]
$|\mathcal{NSL}_6|=\sum_{i=1}^{22}|\mathfrak{A}_i^6|-\sum_{i=1}|\mathfrak{B}_i^6|$
and
\item[(iii)] $|\mathcal{NSL}_7|=\sum_{i=1}^{564}|\mathfrak{A}_i^7|-\sum_{i=1}|\mathfrak{B}_i^7|$.
\end{description}
\end{enumerate}
\end{mycor}
{\bf Proof}\\
An S-loop is an S-groupoid. Thus by Theorem~\ref{1:3}, we have
S-isomorphy classes and S-isotopy classes. Recall that
$|\mathcal{L}_n|=|\mathcal{SL}_n|+|\mathcal{NSL}_n|-|\mathcal{SL}_n\bigcap\mathcal{NSL}_n|$
but $\mathcal{SL}_n\bigcap\mathcal{NSL}_n=\emptyset$ so
$|\mathcal{L}_n|=|\mathcal{SL}_n|+|\mathcal{NSL}_n|$. As stated and
shown in \cite{phd3}, \cite{phd42}, \cite{phd155} and \cite{phd156},
the facts in Table~\ref{table1} are true where $n$ is the order of a
finite loop.
\begin{table}[!hbp]
\begin{center}
\begin{tabular}{|c||c|c|c|}
\hline
$n$ & 5 & 6 & 7  \\
\hline
$|\mathcal{L}_n|$ & 56 & 9, 408 & 16, 942, 080  \\
\hline
$\{\mathcal{A}_i^n\}_{i=1}^k$ & $k=6$ & $k=109$ & $k=23, 746$ \\
\hline
$\{\mathfrak{A}_i^n\}_{i=1}^m$ & $m=2$ & $m=22$ & $m=564$  \\
\hline
\end{tabular}
\end{center}
\caption{Enumeration of Isomorphy and Isotopy classes of finite
loops of small order}\label{table1}
\end{table}
Hence the claims follow.

\begin{myquest}
How many S-loops are in the family $\mathcal{L}_n$? That is, what is
$|\mathcal{SL}_n|$ or $|\mathcal{NSL}_n|$.
\end{myquest}

\begin{myth}\label{1:4.1}
Let $(G,\cdot )$ be a finite S-groupoid of order $n$ with a finite
S-subsemigroup $(H,\cdot )$ of order $m$. Also, let
\begin{displaymath}
\mathcal{ISOT}(G,\cdot ),~\mathcal{SISOT}(G,\cdot
)~\textrm{and}~\mathcal{NSISOT}(G,\cdot )
\end{displaymath}
be the sets of all isotopisms, S-isotopisms and non S-isotopisms of
$(G,\cdot )$. Then,
\begin{displaymath}
\mathcal{ISOT}(G,\cdot )~\textrm{is a group
and}~\mathcal{SISOT}(G,\cdot )\le \mathcal{ISOT}(G,\cdot ).
\end{displaymath}
Furthermore:
\begin{enumerate}
\item $|\mathcal{ISOT}(G,\cdot )|=(n!)^3$;
\item $|\mathcal{SISOT}(G,\cdot )|=(m!)^3$;
\item $|\mathcal{NSISOT}(G,\cdot )|=(n!)^3-(m!)^3$.
\end{enumerate}
\end{myth}
{\bf Proof}\\
\begin{enumerate}
\item This has been shown to be true in [Theorem~4.1.1,
\cite{phd49}].
\item An S-isotopism is an isotopism. So, $\mathcal{SISOT}(G,\cdot
)\subset\mathcal{ISOT}(G,\cdot )$. Thus, we need to just verify the
axioms of a group to show that $\mathcal{SISOT}(G,\cdot
)\le\mathcal{ISOT}(G,\cdot )$. These can be done using the proofs of
reflexivity, symmetry and transitivity in Theorem~\ref{1:3} as
guides. For all triples
\begin{displaymath}
\alpha\in\mathcal{SISOT}(G,\cdot )~\textrm{such that}~\alpha
=(U,V,W)~:~(G,\cdot )\longrightarrow (G,\circ ),
\end{displaymath}
where $(G,\cdot )$ and $(G,\circ )$ are S-groupoids with S-subgroups
$(H,\cdot )$ and $(K,\circ )$ respectively, we can set
\begin{displaymath}
U':=U|_H,~V':=V|_H~\textrm{and}~W':=W|_H~\textrm{since}~A(H)=K~\forall~A\in
\{U,V,W\},
\end{displaymath}
so that $\mathcal{SISOT}(H,\cdot )=\{(U',V',W')\}$. This is possible
because of the following arguments.

Let
\begin{displaymath}
X=\Big\{f':=f|_H~\Big|~f:G\longrightarrow G,~f:H\longrightarrow
K~\textrm{is bijective and}~f(H)=K\Big\}.
\end{displaymath}
Let
\begin{displaymath}
SYM(H,K)=\{\textrm{bijections from $H$ unto $K$}\}.
\end{displaymath}
By definition, it is easy to see that $X\subseteq SYM(H,K)$. Now,
for all $U\in SYM(H,K)$, define $U~:~H^c\longrightarrow K^c$ so that
$U~:~G\longrightarrow G$ is a bijection since $|H|=|K|$ implies
$|H^c|=|K^c|$. Thus, $SYM(H,K)\subseteq X$ so that $SYM(H,K)=X$.

Given that $|H|=m$, then it follows from (1) that
\begin{displaymath}
|\mathcal{ISOT}(H,\cdot )|=(m!)^3~\textrm{so that}~
|\mathcal{SISOT}(G,\cdot )|=(m!)^3~\textrm{since}~SYM(H,K)=X.
\end{displaymath}
\item
\begin{displaymath}
\mathcal{NSISOT}(G,\cdot )=\big(\mathcal{SISOT}(G,\cdot )\big)^c.
\end{displaymath}
So, the identity isotopism
\begin{displaymath}
(I,I,I)\not\in\mathcal{NSISOT}(G,\cdot ),~\textrm{hence}~
\mathcal{NSISOT}(G,\cdot )\nleq\mathcal{ISOT}(G,\cdot ).
\end{displaymath}
Furthermore,
\begin{displaymath}
|\mathcal{NSISOT}(G,\cdot )|=(n!)^3-(m!)^3.
\end{displaymath}
\end{enumerate}

\begin{mycor}\label{1:4.2}
Let $(G,\cdot )$ be a finite S-groupoid of order $n$ with an
S-subsemigroup $(H,\cdot )$. If $\mathcal{ISOT}(G,\cdot )$ is the
group of all isotopisms of $(G,\cdot )$ and $S_n$ is the symmetric
group of degree $n$, then
\begin{displaymath}
\mathcal{ISOT}(G,\cdot )\succsim S_n\times S_n\times S_n.
\end{displaymath}
\end{mycor}
{\bf Proof}\\
As concluded in [Corollary~1, \cite{phd49}], $\mathcal{ISOT}(G,\cdot
)\cong S_n\times S_n\times S_n$. Let $\mathcal{PISOT}(G,\cdot )$ be
the set of all principal isotopisms on $(G,\cdot )$.
$\mathcal{PISOT}(G,\cdot )$ is an S-subgroup in
$\mathcal{ISOT}(G,\cdot )$ while $S_n\times S_n\times \{I\}$ is an
S-subgroup in $S_n\times S_n\times S_n$. If
\begin{displaymath}
\Upsilon~:~\mathcal{ISOT}(G,\cdot )\longrightarrow S_n\times
S_n\times S_n~\textrm{is defined
as}
\end{displaymath}
\begin{displaymath}
\Upsilon\big((A,B,I)\big)=<A,B,I>~\forall~(A,B,I)\in\mathcal{ISOT}(G,\cdot
),
\end{displaymath}
then
\begin{displaymath}
\Upsilon\Big(\mathcal{PISOT}(G,\cdot )\Big)=S_n\times S_n\times
\{I\}.~\therefore~\mathcal{ISOT}(G,\cdot )\succsim S_n\times
S_n\times S_n.
\end{displaymath}

\section{Smarandache $f,g$-Isotopes of Smarandache Loops}
\begin{myth}\label{1:7}
Let $(G,\cdot )$ and $(H,\ast )$ be S-groupoids. If $(G,\cdot )$ and
$(H,\ast )$ are S-isotopic, then $(H,\ast )$ is S-isomorphic to some
Smarandache principal isotope $(G,\circ )$ of $(G,\cdot )$.
\end{myth}
{\bf Proof}\\
Since $(G,\cdot )$ and $(H,\ast )$ are S-isotopic S-groupoids with
S-subsemigroups $(G_1,\cdot )$ and $(H_1,\ast )$, then there exist
bijections $U, V, W : (G,\cdot )\rightarrow (H,\ast )$ such that the
triple $\alpha =(U, V, W)~:~(G,\cdot )\rightarrow (H,\ast )$ is an
isotopism and $\big (G_1\big) A=H_1~\forall~A\in\{U,V,W\}$. To prove
the claim of this theorem, it suffices to produce a closed binary
operation '$\ast$' on $G$, bijections $X,Y~:~G\rightarrow G$, and
bijection $Z~:~G\rightarrow H$ so that
\begin{itemize}
\item the triple $\beta =(X,Y,I)~:~(G,\cdot )\rightarrow (G,\circ )$
is a Smarandache principal isotopism and
\item $Z~:~(G,\circ )\rightarrow (H,\ast )$ is an S-isomorphism or
the triple $\gamma =(Z,Z,Z)~:~(G,\circ )\rightarrow (H,\ast )$ is an
S-isotopism.
\end{itemize}
Thus, we need $(G,\circ )$ so that the commutative diagram below is
true:
\begin{diagram}
(G,\cdot ) &\rTo^{\alpha}_{\textrm{isotopism}} &(H,\ast ) \\
&\rdTo^{\beta}_{\textrm{principal isotopism}} &\uTo^{\gamma}_{\textrm{isomorphism}}\\
& &(G,\circ)
\end{diagram}
because following the proof of transitivity in Theorem~\ref{1:3},
$\alpha =\beta\gamma$ which implies $(U, V, W)=(XZ,YZ,Z)$ and so we
can make the choices; $Z=W$, $Y=VW^{-1}$, and $X=UW^{-1}$ and
consequently,
\begin{displaymath}
x\cdot y=xUW^{-1}\circ VW^{-1}~\Longleftrightarrow~x\circ
y=xWU^{-1}\cdot yWV^{-1}~\forall~x,y\in G.
\end{displaymath}
Hence, $(G,\circ )$ is a groupoid principal isotope of $(G,\cdot )$
and $(H,\ast )$ is an isomorph of  $(G,\circ )$. It remains to show
that these two relationships are Smarandache.

Note that $\big((H_1)Z^{-1}, \circ\big)=(G_1,\circ )$ is a non-trivial
subsemigroup in $(G,\circ )$. Thus, $(G,\circ )$ is an S-groupoid.
So $(G,\circ )\succsim (H,\ast )$. $(G,\cdot )$ and $(G,\circ )$ are Smarandache principal isotopes
because $(G_1)UW^{-1}=(H_1)W^{-1}=(H_1)Z^{-1}=G_1$ and
$(G_1)VW^{-1}=(H_1)W^{-1}=(H_1)Z^{-1}=G_1$.

\begin{mycor}\label{1:8} Let $(G,\cdot )$ be an S-groupoid with an arbitrary groupoid isotope $(H,\ast
)$. Any such groupoid $(H,\ast )$ is an S-groupoid if and only if
all the principal isotopes of $(G,\cdot )$ are S-groupoids.
\end{mycor}
{\bf Proof}\\
By classical result in principal isotopy [\cite{phd3},
III.1.4~Theorem], if $(G,\cdot )$ and $(H,\ast )$ are isotopic
groupoids, then $(H,\ast )$ is isomorphic to some principal isotope
$(G,\circ )$ of $(G,\cdot )$. Assuming $(H,\ast )$ is an S-groupoid
then since $(H,\ast )\cong (G,\circ )$, $(G,\circ )$ is an
S-groupoid. Conversely, let us assume all the principal isotopes of
$(G,\cdot )$ are S-groupoids. Since $(H,\ast )\cong (G,\circ )$,
then $(H,\ast )$ is an S-groupoid.

\begin{myth}\label{1:9}
Let $(G,\cdot )$ be an S-quasigroup. If $(H,\ast )$ is an S-loop
which is S-isotopic to $(G,\cdot )$, then there exist S-elements $f$
and $g$ so that $(H,\ast )$ is S-isomorphic to a Smarandache $f,g$
principal isotope $(G,\circ )$ of $(G,\cdot )$.
\end{myth}
{\bf Proof}\\
An S-quasigroup and an S-loop are S-groupoids. So by
Theorem~\ref{1:7}, $(H,\ast )$ is S-isomorphic to a Smarandache
principal isotope $(G,\circ )$ of $(G,\cdot )$. Let $\alpha
=(U,V,I)$ be the Smarandache principal isotopism of $(G,\cdot )$
onto $(G,\circ )$. Since $(H,\ast )$ is a S-loop and $(G,\circ
)\succsim (H,\ast )$ implies that $(G,\circ )\cong (H,\ast )$, then
$(G,\circ )$ is necessarily an S-loop and consequently, $(G,\circ )$
has a two-sided identity element say $e$ and an S-subgroup
$(G_2,\circ )$. Let $\alpha =(U,V,I)$ be the Smarandache principal
isotopism of $(G,\cdot )$ onto $(G,\circ )$. Then,
\begin{displaymath}
xU\circ yV=x\cdot y~\forall~x,y\in G\Longleftrightarrow~x\circ
y=xU^{-1}\cdot yV^{-1}~\forall~x,y\in G.
\end{displaymath}
So,
\begin{displaymath}
y=e\circ y=eU^{-1}\cdot yV^{-1}=yV^{-1}L_{eU^{-1}}~\forall~y\in
G~\textrm{and}~x=x\circ e=xU^{-1}\cdot
eV^{-1}=xU^{-1}R_{eV^{-1}}~\forall~x\in G.
\end{displaymath}
Assign $f=eU^{-1},g=eV^{-1}\in G_2$. This assignments are well
defined and hence $V=L_f$ and $U=R_g$. So that $\alpha =(R_g,L_f,I)$
is a Smarandache $f,g$ principal isotopism of $(G,\circ )$ onto
$(G,\cdot )$. This completes the proof.

\begin{mycor}\label{1:10}
Let $(G,\cdot )$ be an S-quasigroup(S-loop) with an arbitrary
groupoid isotope $(H,\ast )$. Any such groupoid $(H,\ast )$ is an
S-quasigroup(S-loop) if and only if all the principal isotopes of
$(G,\cdot )$ are S-quasigroups(S-loops).
\end{mycor}
{\bf Proof}\\
This follows immediately from Corollary~\ref{1:8} since an
S-quasigroup and an S-loop are S-groupoids.

\begin{mycor}\label{1:11}
If $(G,\cdot )$ and  $(H,\ast )$ are S-loops which are S-isotopic,
then there exist S-elements $f$ and $g$ so that $(H,\ast )$ is
S-isomorphic to a Smarandache $f,g$ principal isotope $(G,\circ )$
of $(G,\cdot )$.
\end{mycor}
{\bf Proof}\\
An S-loop is an S-quasigroup. So the claim follows from
Theorem~\ref{1:9}.

\section{G-Smarandache Loops}
\begin{mylem}\label{1:12}
Let $(G,\cdot )$ and $(H,\ast )$ be S-isotopic S-loops. If $(G,\cdot
)$ is a group, then $(G,\cdot )$ and $(H,\ast )$ are S-isomorphic
groups.
\end{mylem}
{\bf Proof}\\
By Corollary~\ref{1:11}, there exist S-elements $f$ and $g$ in
$(G,\cdot )$ so that $(H,\ast )\succsim (G,\circ )$ such that
$(G,\circ )$ is a Smarandache $f,g$ principal isotope of $(G,\cdot
)$.

Let us set the mapping $\psi :=R_{f\cdot g}=R_{fg}~:~G\rightarrow
G$. This mapping is bijective. Now, let us consider when $\psi
:=R_{fg}~:~(G,\cdot )\rightarrow (G,\circ )$. Since $(G,\cdot )$ is
associative and $x\circ y=xR_g^{-1}\cdot yL_f^{-1}~\forall~x,y\in
G$, the following arguments are true.

$x\psi \circ y\psi =x\psi R_g^{-1}\cdot y\psi L_f^{-1}=xR_{fg}
R_g^{-1}\cdot yR_{fg}L_f^{-1}=x\cdot fg\cdot g^{-1}\cdot f^{-1}\cdot
y\cdot fg=x\cdot y\cdot fg=(x\cdot y)R_{fg}=(x\cdot
y)\psi~\forall~x,y\in G$. So, $(G,\cdot )\cong (G,\circ )$. Thus,
$(G,\circ )$ is a group. If $(G_1,\cdot )$ and $(G_1,\circ )$ are
the S-subgroups in $(G,\cdot )$ and $(G,\circ )$, then
$\big((G_1,\cdot )\big) R_{fg}=(G_1,\circ )$. Hence, $(G,\cdot
)\succsim (G,\circ )$.

$\therefore~(G,\cdot )\succsim (H,\ast )$ and $(H,\ast )$ is a
group.

\begin{mycor}\label{1:13}
Every group which is an S-loop is a GS-loop.
\end{mycor}
{\bf Proof}\\
This follows immediately from Lemma~\ref{1:12} and the fact that a
group is a G-loop.

\begin{mycor}\label{1:14}
An S-loop is S-isomorphic to all its S-loop S-isotopes if and only
if it is S-isomorphic to all its Smarandache $f,g$ principal
isotopes.
\end{mycor}
{\bf Proof}\\
Let $(G,\cdot )$ be an S-loop with arbitrary S-isotope $(H,\ast )$.
Let us assume that $(G,\cdot )\succsim (H,\ast )$. From
Corollary~\ref{1:11}, for any arbitrary S-isotope $(H,\ast )$ of
$(G,\cdot )$, there exists a Smarandache $f,g$ principal isotope
$(G,\circ )$ of $(G,\cdot )$ such that $(H,\ast )\succsim (G,\circ
)$. So, $(G,\cdot )\succsim (G,\circ )$.

Conversely, let $(G,\cdot )\succsim (G,\circ )$, using the fact in
Corollary~\ref{1:11} again, for any arbitrary S-isotope $(H,\ast )$
of $(G,\cdot )$, there exists a Smarandache $f,g$ principal isotope
$(G,\circ )$ of $(G,\cdot )$ such that $(G,\circ )\succsim (H,\ast
)$. Therefore, $(G,\cdot )\succsim (H,\ast )$.

\begin{mycor}\label{1:15}
A S-loop is a GS-loop if and only if it is S-isomorphic to all its
Smarandache $f,g$ principal isotopes.
\end{mycor}
{\bf Proof}\\
This follows by the definition of a GS-loop and
Corollary~\ref{1:14}.


\begin{thebibliography}{99}
\bibitem{phd41} R. H. Bruck (1966), {\it A survey of binary systems}, Springer-Verlag, Berlin-G\"ottingen-Heidelberg, 185pp.
\bibitem{phd155} R. E. Cawagas (2000), {\it Generation of NAFIL
loops of small order}, Quasigroups and Related Systems, 7, 1--5.
\bibitem{phd39} O. Chein, H. O. Pflugfelder and J. D. H. Smith (1990), {\it Quasigroups and loops : Theory and applications}, Heldermann Verlag, 568pp.
\bibitem{phd49} J. D\'{e}ne and A. D. Keedwell (1974), {\it Latin squares and their applications}, the English University press Lts, 549pp.
\bibitem{phd42} E. G. Goodaire, E. Jespers and C. P. Milies (1996), {\it Alternative loop rings}, NHMS(184), Elsevier, 387pp.
\bibitem{sma1} T. G. Ja\'iy\'e\d ol\'a (2006), {\it An holomorphic study of the Smarandache concept in
loops}, Scientia Magna Journal, 2, 1, 1--8.
\bibitem{sma2} T. G. Ja\'iy\'e\d ol\'a (2006), {\it Parastrophic invariance of Smarandache quasigroups}, Scientia Magna Journal, 2, 3, 48--53.
\bibitem{sma3} T. G. Ja\'iy\'e\d ol\'a (2006), {\it On the universality of some Smarandache loops of Bol-Moufang type}, Scientia Magna Journal, 2, 4, 45--48.
\bibitem{phd156} B. D. McKay, A. Meynert and W. Myrvold (2007), {\it Small latin squares, quasigroups and
loops}, Journal of Combinatorial Designs, 15, 2, 98--119.
\bibitem{muk} A. S. Muktibodh (2006), {\it Smarandache quasigroups},
Scientia Magna Journal, 2, 1, 13--19.
\bibitem{phd3} H. O. Pflugfelder (1990), {\it Quasigroups and loops : Introduction}, Sigma series in Pure Math. 7, Heldermann Verlag, Berlin, 147pp.
\bibitem{phd75} W. B. Vasantha Kandasamy (2002), {\it Smarandache
loops}, Department of Mathematics, Indian Institute of Technology,
Madras, India, 128pp.
\bibitem{phd83} W. B. Vasantha Kandasamy (2002), {\it Smarandache
Loops}, Smarandache Notions Journal, 13, 252--258.
\bibitem{van1} W. B. Vasantha Kandasamy (2002), {\it Groupoids and Smarandache Groupoids},
American Research Press Rehoboth, 114pp.
\bibitem{van2} W. B. Vasantha Kandasamy (2002), {\it Smarandache Semigroups},
American Research Press Rehoboth, 94pp.
\bibitem{van3} W. B. Vasantha Kandasamy (2002), {\it Smarandache Semirings, Semifields, And Semivector Spaces},
American Research Press Rehoboth, 121pp.
\bibitem{van4} W. B. Vasantha Kandasamy (2003), {\it Linear Algebra And Smarandache Linear Algebra},
American Research Press, 174pp.
\bibitem{van5} W. B. Vasantha Kandasamy (2003), {\it Bialgebraic Structures And Smarandache Bialgebraic Structures},
American Research Press Rehoboth, 271pp.
\bibitem{van6} W. B. Vasantha Kandasamy and F. Smarandache (2005), {\it N-Algebraic Structures And Smarandache N-Algebraic Structures},
Hexis Phoenix, Arizona, 174pp.
\end{thebibliography}
\end{document}